 \documentclass[preprint,review,12pt]{elsarticle}

\usepackage{amssymb}

\usepackage{amsmath}
\usepackage{graphicx}
\usepackage{graphics}
\usepackage{subfigure}
\usepackage{mathptmx}
\usepackage[misc]{ifsym}
\usepackage{algorithm}
\usepackage{algorithmicx}
\usepackage{algpseudocode}
\usepackage{multirow}
\usepackage{amsmath}
\usepackage{amsfonts}
\usepackage{graphics}
\usepackage{epsfig}
\usepackage{amssymb}
\usepackage{color}

\usepackage{amsmath}
\usepackage{graphicx}
\usepackage{graphics}
\usepackage{subfigure}
\usepackage{mathptmx}
\usepackage[misc]{ifsym}
\usepackage{algorithm}
\usepackage{algorithmicx}
\usepackage{algpseudocode}
\usepackage{graphicx}
\usepackage{graphics}
\usepackage{subfigure}
\usepackage{epsfig}
\usepackage{color}
\usepackage{multirow}

\usepackage[top=2cm, bottom=2cm, left=3cm, right=3cm] {geometry}

\newtheorem{theorem}{Theorem}
\newtheorem{lemma}{Lemma}

\newproof{proof}{Proof}

\newtheorem{assumption}{Assumption}

\begin{document}

\begin{frontmatter}



\title{{\bf Distributed Regularized Dual Gradient Algorithm for Constrained Convex Optimization over Time-Varying Directed Graphs}}


\author[cq]{Chuanye Gu}  \ead{guchuanye11@163.com}
\author[cq]{Zhiyou Wu} \ead{zywu@cqnu.edu.cn}
\author[cq]{Jueyou Li
} \ead{lijueyou@163.com}
\address[cq]{School of Mathematical Sciences, Chongqing Normal University, Chongqing, 400047, China}

\begin{abstract}

We investigate a distributed optimization problem over a cooperative multi-agent time-varying network, where each agent has its own decision variables that should be set so as to minimize its individual objective subject to local constraints and global coupling constraints.
Based on push-sum protocol and dual decomposition, we design a distributed regularized dual gradient algorithm to solve this problem, in which the algorithm is implemented in time-varying directed graphs only requiring the column stochasticity of  communication matrices.
By augmenting the corresponding Lagrangian function with a quadratic regularization term, we first obtain the bound of the Lagrangian multipliers which does not require constructing a compact set containing the dual optimal set when compared with most of primal-dual based methods.
Then, we obtain that the convergence rate of the proposed method can achieve the order of $\mathcal{O}(\ln T/T)$ for strongly convex objective functions, where $T$ is the iterations. Moreover, the explicit bound of constraint violations is also given.
Finally, numerical results on the network utility maximum problem are used to demonstrate the efficiency of the proposed algorithm.

\end{abstract}

\begin{keyword}
 Convex optimization \sep Distributed algorithm \sep Dual decomposition \sep Regularization \sep Multi-agent network.

\end{keyword}
\end{frontmatter}


\section{Introduction}\label{1}

In recent years it is witnessed the unprecedented growth in the research for solving many optimization problems over multi-agent networks \cite{nedic2009distributed,nedic2010constaints,jakovetic2014fast,nedic2015distributed}. Distributed optimization has been found in a lot of application domains,
such as distributed finite-time optimal rendezvous problem \cite{johansson2008subgradient}, wireless and social networks \cite{baingana2014proximalgradient}, \cite{mateos2012distributed}, power systems \cite{bolognani2015distributed}, \cite{zhangdistributed2016}, robotics \cite{martinea2007on}, and so on. There is indeed a long history in the optimization community of this problem, see \cite{tsitsiklis1986distributed}.

Based on consensus schemes, there are mainly three categories of algorithms designed for distributed optimization in the literatures, including primal consensus distributed algorithms, dual consensus distributed algorithms and primal-dual consensus distributed algorithms, see \cite{nedic2009distributed,ram2010distributed,c2012dual,zhuon2012,li2015gradient,lorenzo2016next}.
In most to previous works, the communication graphs are required to be balanced, i.e., the communication weight matrices are doubly stochastic. The paper \cite{gharesifard2012distributed}
 considered a fixed and directed graph with the requirement of a balanced graph.
The work in \cite{iconsensus2012} proposed distributed subgradient based algorithms in directed and fixed topologies, in which the messages among agents are propagated by ``push-sum" protocol. However, the communication protocol is required to know the number of agents or the graph.
In general, push-sum protocol is attractive for implementations since it can easily operate over directed communication topologies, and thus avoids incidents of deadlock that may occur in practice when using undirected communication topologies \cite{nedic2015distributed}.
Nedi\'{c} et al. in \cite{nedic2015distributed} designed subgradient-push distributed method for a class of unconstrained optimization problems, in which the requirement of a balanced graph was canceled. Their proposed method has a slower convergence rate with order of $\mathcal{O}(\ln T/\sqrt{T})$. Later, Nedi\'{c} et al. in \cite{nedic2015stochastic} improved the convergence rate from $\mathcal{O}(\ln T/\sqrt{T})$ to $\mathcal{O}(\ln T/T)$ under the condition of strong convexity. However, they only considered unconstrained optimization problems.

The methods for solving distributed optimization problems subject to equality or (and) inequality constraints have received considerable attention \cite{bertsekas2003convex,necoara2008application,li2016a}. The authors in \cite{zhuon2012} first proposed a distributed Lagrangian primal-dual subgradient method by characterizing the primal-dual optimal solutions as the saddle points of the Lagrangian function related to the problem under consideration. The work \cite {yuan2011distributed} developed a variant of the distributed primal-dual subgradient method by introducing multistep consensus mechanism. For more general distributed optimization problem with inequality constraints that couple all the agents' decision variables, Chang et al. \cite{chang2014distributed} designed a novel distributed primal-dual perturbed subgradient method and analyzed the convergence. The implementation of the algorithms aforementioned usually involves projections onto some primal and dual constrained sets, respectively. In particular, they require constructing a compact set that contains the dual optimal set, and projecting the dual variable onto this set to guarantee the boundedness of dual iterates, which is of importance in establishing the convergence of the algorithms. However, the construction of this compact set is impractical since it involves each agent solving a general constrained convex problem \cite{Yuan2016Regularized,Khuzani2016Distributed}.
To ensure the boundedness of the norm of the dual variables, Yuan et al. in \cite{Yuan2016Regularized} proposed a regularized primal-dual distributed algorithm. However, the optimization problem only includes one constraint. Later, Khuzani et al. in \cite{Khuzani2016Distributed} investigated a distributed optimization with several inequality constraints, and established the convergence of their proposed distributed deterministic and stochastic primal-dual algorithms, respectively. Very recently, Falsone et al. \cite{falsone2016dual} designed a dual decomposition based distributed method for solving a separable convex optimization with coupled inequality constraints and provided the convergence analysis, but none of explicit convergence rate of their algorithm was given. Most of aforementioned works operating over undirected networks with the usage of doubly stochastic matrices are possible. However, it turns out that directed graphs depending on doubly stochastic matrices may be undesirable for a variety of reasons, see \cite{nedic2015distributed,nedic2015stochastic}.

In this paper, we propose a distributed regularized dual gradient method for solving convex optimization problem subjected to local and coupling constraints over time-varying directed networks.
The proposed method is based on push-sum protocol. Each agent is only required to know its out-degree at each time, without requiring knowledge of either the number of agents or the graph sequence. By augmenting the corresponding Lagrangian function with a quadratic regularization term,  the norm of the multipliers is bounded, which does not require constructing a compact set containing
the dual optimal set when compared with existing most of primal-dual methods. The convergence rate of the method with the order of $\mathcal{O}(\ln t/t)$ for strongly convex objective functions is obtained. Moreover, the explicit bound on the constraint violations is also provided.

The main contributions of this paper are two folds. Firstly, we establish the upper bound on the norm of dual variables by resorting to the regularized Lagrangian function. Secondly, we obtain the explicit convergence rates of the proposed method over the directed unbalanced network.
The work in this paper is related to the recent literatures  \cite{nedic2015stochastic}  and \cite{falsone2016dual}. The reference in \cite{nedic2015stochastic} addresses an unconstrained distributed optimization over time-varying directed networks, while our paper investigates a distributed optimization with coupling equality constraints. Our method can be viewed as an extension of push-sum based algorithms \cite{nedic2015stochastic} to a constrained setting.
Compared with the method in  \cite{falsone2016dual}, our proposed distributed algorithm is inspired by push-sum strategy over time-varying directed networks without the requirement of balanced network graphs, whereas the method in \cite{falsone2016dual} must require that the graphs are balanced and the communication matrices  are doubly stochastic. In \cite{falsone2016dual}, the authors only establish the convergence of their approach. However, in this paper, we obtain the explicit convergence rates of the proposed method in the time-varying directed network topology. More importantly, we further give the explicit convergence estimate on constraint violations.
The regularized primal-dual distributed methods proposed in \cite{Yuan2016Regularized,Khuzani2016Distributed} require that the networks are undirected and the communication weight matrices are double stochastic, whereas our method can deal with distributed optimization problems over time-varying directed graphs, only needing the column stochastic matrices.

The remainder of this paper is organized as follows. In Section \ref{2}, we state the related problem, useful assumptions and preparatory work.
In Section \ref{3}, we propose the distributed regularized dual gradient algorithm and give main results. In Section \ref{4}, we give some Lemmas and the proof of main results. Numerical simulations are given in Section \ref{5}. Finally, Section \ref{6} draws some conclusions.

Notation: We use boldface to distinguish between the scalars and vectors in $\mathbb{R}^{n}$.
For example, $v_{i}[t]$ is a scalar and $\mathbf{u}_{i}[t]$ is a vector. For a matrix $W$, we will use the $(W)_{ij}$ to show its $i,j$'th entry.
We use the $||\mathbf{x}||$ to denote the Euclidean norm of a vector $\mathbf{x}$, and $\mathbf{1}$ for the vector of ones.
A convex function $f: \mathbb{R}^{n}\rightarrow \mathbb{R}$ is $\widetilde{\gamma}$-strongly convex with $\widetilde{\gamma}>0$ if the following relation holds, for all $\mathbf{x}, \mathbf{y}\in  \mathbb{R}^{n}$
$$f(\mathbf{x})-f(\mathbf{y}) \geq \nabla g(\mathbf{y})^{\top}(\mathbf{x}-\mathbf{y}) + \frac{\widetilde{\gamma}}{2}||\mathbf{x}-\mathbf{y}||^2,$$
where $g(\mathbf{y})$ is any subgradient of $f$ at $\mathbf{y}$.

\section{Distributed optimization problem with equality constraints} \label{2}

\subsection{Constrained Multi-agent Optimization }

Consider the following constrained optimization problem
\begin{equation}\label{cop}
   \min_{\{\mathbf{x}_{i}\in \mathbf{X}_{i}\}_{i=1}^{m}} ~~F(\mathbf{x}):=\sum_{i=1}^{m} f_{i}(\mathbf{x}_{i}) ~~~~\textrm{s.t.}~~~~ \sum_{i=1}^{m}(A_{i}\mathbf{x}_{i}-\mathbf{b}_{i})=0,
\end{equation}
where there are $m$ agents associated with a time-varying network.
Each agent $i$ only knows its own objective function $f_{i}(\mathbf{x}_{i})$: $\mathbb{R}^{n_{i}}\rightarrow \mathbb{R}$ and its own constraints $\mathbf{X}_{i}\in \mathbb{R}^{n_{i}}$, and all agents subject to the coupling equality constraints $\sum_{i=1}^{m}(A_{i}\mathbf{x}_{i}-\mathbf{b}_{i})=0$, $A_{i}\in \mathbb{R}^{p\times n_i}$ and $\mathbf{b}_{i}\in \mathbb{R}^{p}$. $\mathbf{x}=(\mathbf{x}_{1}^{\top},\mathbf{x}_{2}^{\top},\cdots, \mathbf{x}_{m}^{\top})^{\top}$ with $n=\sum_{i=1}^{m} n_{i}$, belongs to $\mathbf{X}=\mathbf{X}_{1}\times \mathbf{X}_{2}\times\cdots\times \mathbf{X}_{m}$.

Problem (\ref{cop}) is quite general arising in diverse applications, for examples, distributed model predictive control \cite{necoara2008application}, network utility maximization \cite{Low1999optimization,beck2014an}, economic dispatch problems for smart grid \cite{bolognani2015distributed,zhangdistributed2016}.

To decouple the coupling equality constraints, we introduce a regularized Lagrangian function $\mathcal{L}(\mathbf{x},\lambda)$ of problem (\ref{cop}), given by
\begin{eqnarray}\label{eq-regularized-L-fuction}
\mathcal{L} (\mathbf{x},\lambda) :=
\sum_{i=1}^{m} [f_{i}(\mathbf{x}_{i}) + \lambda^{\top}(A_{i}\mathbf{x}_{i}-\mathbf{b}_{i})
 - \frac{\gamma_{i}}{2}\lambda^{\top}\lambda]=\sum_{i=1}^{m} \mathcal{L}_{i}(\mathbf{x}_{i},\lambda) ,
\end{eqnarray}
where $\mathcal{L}_{i}(\mathbf{x}_{i},\lambda)=f_{i}(\mathbf{x}_{i}) + \lambda^{\top}(A_{i}\mathbf{x}_{i}-\mathbf{b}_{i})
 - \frac{\gamma_{i}}{2}\lambda^{\top}\lambda$ are the regularized Lagrangian function associated with the $i$th agent, and $\gamma_{i}>0$ is regularization parameter, for $i=1,2,\ldots,m$.

Define a regularized dual function of problem (\ref{cop}) as follows
$$\phi(\lambda):= \min_{\mathbf{x}\in \mathbf{X}}\mathcal{L}(\mathbf{x},\lambda).$$
Note that the regularized Lagrangian function $\mathcal{L}(\mathbf{x},\lambda)$ defined in (\ref{eq-regularized-L-fuction}) is separable with respect to $\mathbf{x}_{i}, i=1,\ldots,m$. Thus, the regularized dual function $\phi(\lambda)$ can be rewritten as
\begin{equation}\label{li}
  \phi(\lambda)= \sum_{i=1}^{m} \phi_{i}(\lambda)
  =\sum_{i=1}^{m}\min_{\mathbf{x}_{i}\in \mathbf{X}_{i}}\mathcal{L}_{i}(\mathbf{x}_{i},\lambda),
\end{equation}
where $\phi_{i}(\lambda):=\min_{\mathbf{x}_{i}\in \mathbf{X}_{i}}\mathcal{L}_{i}(\mathbf{x}_{i},\lambda)$ can be regarded as the regularized dual function of agent $i, i=1,\ldots,m$.

Then, the regularized dual problem of problem (\ref{cop}) can be written as $\max_{\lambda}\min_{\mathbf{x}\in X} \mathcal{L}(\mathbf{x},\lambda)$, or, equivalently,
\begin{equation}\label{dual}
 \max_{\lambda}\sum_{i=1}^{m} \phi_{i}(\lambda).
\end{equation}

The coupling equality constraints between agents is represented by the fact that $\lambda$ is a common decision vector and all the agents should agree on its value.

\subsection{Related assumptions }

The following assumptions on the problem (\ref{cop}) and on the communication time-varying network are needed to show properties of convergence for the proposed method.

\begin{assumption}\label{A1}
For each $i=1,2,\ldots, m$, the function $f_{i}(\cdot)$: $\mathbb{R}^{n_{i}}\rightarrow \mathbb{R}$ is $\tau_{i}$ strongly convex, and the set $\mathbf{X}_{i}\subseteq \mathbb{R}^{n_{i}}$ is non-empty, convex and compact.
\end{assumption}

Note that, under the Assumption \ref{A1}, we have:

(i) the function $\phi_{i}(\lambda)$ defined in (\ref{dual}) is $\gamma_{i}$-strongly concave, differentiable and its gradient
$\nabla\phi_{i}(\lambda)=A_{i}\mathbf{x}_{i}(\lambda)-\mathbf{b}_{i}-\gamma_{i} \lambda$ is Lipschitz continuous with constant $||A_{i}||/\tau_{i}$, where $\mathbf{x}_{i}(\lambda):=\arg\min_{\mathbf{x}_{i}\in \mathbf{X}_{i}} \mathcal{L}_{i}(\mathbf{x}_{i},\lambda)$ (see \cite{beck2014an,li2016a}, for more details);

(ii) for any $\mathbf{x}_{i}\in \mathbf{X}_{i}$, there is a constant $G_{i}>0$ such that $||A_{i}\mathbf{x}_{i}-\mathbf{b}_{i}||\leq G_{i}$, due to the compactness of $\mathbf{X}_{i}$, $i=1,2,\ldots, m$. 

We assume that each agent can communicate with other agents over a time-varying network. The communication topology is modeled by a directed graph $\mathcal{G}[t]=(\mathcal{V}, \mathcal{E}[t])$ over the vertex set $\mathcal{V}=\{1,\ldots,m\}$ with the edge set $\mathcal{E}[t]\subseteq \mathcal{V}\times \mathcal{V}$.
Let $\mathcal{N}_{i}^{in}[t]$ represent the collection of in-neighbors and $\mathcal{N}_{i}^{out}[t]$ represent the collection of out-neighbors of agent $i$ at time $t$, respectively. That is
$$\mathcal{N}_{i}^{in}[t]:=\{j|(j,i)\in \mathcal{E}[t]\}\cup  \{i\},$$
$$\mathcal{N}_{i}^{out}[t]:=\{j|(i,j)\in \mathcal{E}[t]\}\cup  \{i\},$$
where $(j,i)$ represents agent $j$ may send its information to agent $i$.
And let $d_{i}(t)$ be the out-degree of agent $i$, i.e.,
$$d_{i}[t]=|\mathcal{N}_{i}^{out}[t]|,$$

We introduce a time-varying communication weight matrix $W[t]$ with elements $(W[t])_{ij}$, defined by
\begin{eqnarray}\label{cop2}
& (W[t])_{ij}~~ =  ~~\left\{
\begin{aligned}
&~~~~\frac{1}{d_{j}[t]},~~\textrm{when}~ j\in \mathcal{N}_{i}^{in}[t],~i,~j=1,2,\ldots,m,\\
&~~~~~~0,~~~~~~~~\textrm{otherwise}.\\
\end{aligned}
\right.
\end{eqnarray}

We need the following assumption on the weight matrix $W[t]$, which can be found in \cite{nedic2015distributed}, \cite{nedic2010constaints}.
\begin{assumption}\label{A2}
i) Every agent $i$ knows its out-degree $d_{i}[t]$ at every time $t$;~
ii) The graph sequence $\mathcal{G}[t]$ is $B$-strongly connected, namely, there exists an integer $B>0$ such that the sequence $\mathcal{G}[t]$ with edge set $\mathcal{E}[t]=\cup_{l=kB}^{(k+1)B-1}\mathcal{E}[l]$ is strongly connected, for all $t\geq0$.
\end{assumption}

Note that the communicated weight matrix $W[t]$ is column-stochastic. In this paper, we do not require the assumption of double-stochasticity on $W[t]$.

\section{Algorithm and main results }\label{3}

\subsection{Distributed regularized dual gradient algorithm}

In general, the problem (\ref{cop}) could be solved in a centralized manner. However, if the number $m$ of agents is large, this may turn out to be computationally challenge. Additionally, each agent would
be required to share its own information, such as the objective $f_{i}$, the constraints $X_{i}$ and $(A_{i},\mathbf{b}_{i})$, either with the other agents or with a central coordinate collecting
all information, which is possibly undesirable in many cases, due to privacy concerns.

To overcome both the computational and privacy issues stated above, we propose a Distributed Regularized Dual Gradient Algorithm (DRDGA, for short) by resorting to solve the regularized dual problem (\ref{dual}). Our proposed algorithm DRDGA is motivated by the gradient push-sum method  \cite{nedic2015distributed} and dual decomposition \cite{falsone2016dual,li2016a}, described as in Algorithm \ref{alg1}.

\begin{algorithm}[!htb]
\caption{Distributed Regularized Dual Gradient Algorithm (DRDGA) }\label{alg1}
\begin{algorithmic}[1]
\State Initialization: for $i=1,2,\ldots,m$, given $\mathbf{\theta}_{i}[0]\in \mathbb{R}^{p}$, $\rho_{i}[0]=1$, $i=1,2,\ldots,m$; set $t:=0$;
\Repeat
\For{each agent $i=1,\ldots,m$}
\State $\mathbf{u}_{i}[t+1] = \sum_{j=1}^{m}(W[t])_{ij} \mathbf{\theta}_{j}[t]$;
\State $\rho_{i}[t+1] = \sum_{j=1}^{m}(W[t])_{ij} \rho_{j}[t]$;
\State $\lambda_{i}[t+1] = \frac{\mathbf{u}_{i}[t+1]} {\rho_{i}[t+1]}$;
\State $\mathbf{x}_{i}[t+1] =
\arg \min_{\mathbf{x}_{i} \in \mathbf{X}_{i}} \{f_{i}(\mathbf{x}_{i}) +
\lambda_{i}[t+1]^{\top} (A_{i}\mathbf{x}_{i} - \mathbf{b}_{i}) - \frac{\gamma_{i}}{2} \lambda_{i}[t+1]^{\top}\lambda_{i}[t+1]\};$
\State $\mathbf{\theta}_{i}[t+1]=\mathbf{u}_{i}[t+1]+\beta[t+1](A_{i}\mathbf{x}_{i}[t+1]-\mathbf{b}_{i} - \gamma_{i} \lambda_{i}[t+1])$;
\EndFor
\State set $t=t+1$;
\Until{a preset stopping criterion is met.} 
\end{algorithmic}
\end{algorithm}

In Algorithm \ref{alg1}, each agent $i$ broadcasts (or pushes) the quantities $\mathbf{\theta}_{i}[t]/d_{i}[t]$ and $\rho_{i}[t]/d_{i}[t]$ to all of the agents in its out-neighborhood $\mathcal{N}_{i}^{out}[t]$.  Then, each agent simply sums all the received messages to obtain $\mathbf{u}_{i}[t+1]$ in step 4 and $\mathbf{\rho}_{i}[t+1]$ in step 5, respectively. The update rules in steps 6-8 can be implemented locally. In particular, the update of local primal vector $\mathbf{x}_{i}[t+1]$ in step 7 is performed by minimizing $\mathcal{L}_{i}$ with respect to $\mathbf{x}_{i}$ evaluated at $\lambda=\lambda_{i}[t+1]$, while the update of the dual
vector $\mathbf{\lambda}_{i}[t+1]$ in step 8  involves  the
maximization of $\mathcal{L}_{i}$ with respect to $\mathbf{\lambda}=\lambda_{i}$ evaluated at $\mathbf{x}_{i}=\mathbf{x}_{i}[t+1]$.
Note that the term  $A_{i}\mathbf{x}_{i}[t+1]-\mathbf{b}_{i}-\gamma_{i} \lambda_{i}[t+1]$ in step 8 is the gradient of $\phi_{i}(\lambda)$ at $\lambda=\lambda_{i}[t+1]$.

\subsection{Statement of main results}

In this section, we will show that the main results of the convergence for the proposed Algorithm \ref{alg1}.

It is shown in \cite{nedic2010constaints} that the local primal vector $\mathbf{x}_{i}[t]$ does not converge to the optimal solution $\mathbf{x}_{i}^{*}$ of problem (\ref{cop}) in general. Compared to $\mathbf{x}_{i}[t]$, however, the following recursive auxiliary primal iterates
$$\widehat{\mathbf{x}}_{i}[T] = \frac{\sum_{t=1}^{T} (t-1) \mathbf{x}_{i}[t]}{\frac{T(T-1)}{2}}, ~\mathrm{for~all}~T\geq 2$$
can show better convergence properties by setting $\widehat{\mathbf{x}}_{i}[1]=\mathbf{x}_{i}[0]$, see \cite{zhuon2012,chang2014distributed,beck2014an}.
Define the averaging iterates as $\overline{\theta}[t]=\frac{\sum_{i=1}^{m}\mathbf{\theta}_{i}[t]}{m}$.

The following Theorem \ref{th1} first give an upper on the norm of dual variables. By controlling the norm of the dual variables, we in turn control the norm of the sub-gradients of the augmented Lagrangian function, which are instrumental to prove Theorem \ref{th2} and Theorem \ref{th3} below.

\begin{theorem}\label{th1}
Suppose that Assumptions \ref{A1} and \ref{A2} hold and the non-increasing stepsize sequence $\{\beta[t]\}_{t>0}$ satisfies $\lim_{t\rightarrow \infty} \beta[t] = 0$. Then, there is a positive constant $D$ such that for all $i = 1, 2, \ldots, m,$ $$\sup_{t}||\lambda_{i}[t]|| \leq D.$$
\end{theorem}

In what follows, Theorem \ref{th2} shows the convergence rate of primal function's value under  Assumptions \ref{A1} and \ref{A2}.

\begin{theorem}\label{th2} (Convergence rate)
Suppose Assumptions \ref{A1} and \ref{A2} are satisfied and the stepsize is taken as $\beta[t]=\frac{q}{t}, t=1,2,\ldots$, where the constant $q$ is such that $\frac{q\gamma}{m}\geq 4$. Then, for all $T\geq 1$ and $i=1,2,\ldots, m$, we have
\begin{eqnarray}
&&F(\widehat{\mathbf{x}}_{i}[T]) - F(\mathbf{x}^{*}) \nonumber\\
&\leq & \frac{32}{T\delta}\sum_{i=1}^{m}(G_{i} + \gamma_{i}D)\left (\frac{\eta}{ 1 - \eta}\sum_{i=1}^{m}||\theta_{i}[0]||_{1} +
\frac{ q m B}{1-\eta}(1 + \ln T)\right )\nonumber\\
&&+ \frac{q}{T}\sum_{i=1}^{m}(G_{i} + \gamma_{i}D)^{2}.\nonumber
\end{eqnarray}
where $D$ is the bound of dual variable, $B= \max_{1\leq i\leq m} \sqrt{p}(G_{i} + \gamma_{i}D)$, $\gamma=\sum_{j=1}^{m}\gamma_{j}$, the scalar $\eta\in(0,1)$ and $\delta>0$ satisfy $\delta \geq \frac{1}{m^{mB}},~~ \eta \leq(1-\frac{1}{m^{mB}})^{\frac{1}{mB}}.$
\end{theorem}

Theorem \ref{th2} shows that the iterative sequence of primal objective function $\{F(\widehat{\mathbf{x}}[T])\}$ converges to the optimal value $F(\mathbf{x}^{*})$ at a rate of $O(\ln T/T)$, i.e.,
$$F(\widehat{\mathbf{x}}[T])-F(\mathbf{x}^{*})=O\left(\frac{\ln T}{T}\right)$$
with the constant relying on the regularization parameters $\gamma_{i}, i=1,2,\ldots,m$, the bounds of dual variables $D$ and coupling constraints $G_{i}, i=1,2,\ldots,m$, initial values $\overline{\mu}[0]$ at the agents, and on both the speed $\eta$ of the network information diffusion and the imbalance $\delta$ of influences among the agents.

In the next theorem, we show that the upper bound on the constraint violation.

\begin{theorem}\label{th3} (Constraint violation bound)
Suppose Assumptions \ref{A1} and \ref{A2} are satisfied and the stepsize is taken as $\beta[t]=\frac{q}{t}, t=1,2,\ldots$, where the constant $q$ is such that $\frac{q\gamma}{m}\geq 4$. Then, for all $T\geq 1$ and $i=1,2,\ldots, m$, we have
\begin{eqnarray}
&&||\sum_{j=1}^{m}A_{j}\widehat{\mathbf{x}}_{j}[T]-\mathbf{b}_{j}||^{2}\nonumber\\
&& \leq \frac{\gamma}{T\delta}\sum_{j=1}^{m}(G_{j}+ \gamma_{j}D)\left(\frac{8\eta}{1-\eta}\sum_{j=1}^{m}||\mu_{j}[0]||_{1}+\frac{8qmB}{1-\eta}(1+\ln T)\right)\nonumber\\
&& + \frac{q\gamma}{4T}\sum_{j=1}^{m}(G_{j}+ \gamma_{j}D)^{2}.\nonumber
\end{eqnarray}
where $D$ is the bound of dual variable, $B= \max_{1\leq i\leq m} \sqrt{p}(G_{i} + \gamma_{i}D)$, $\gamma=\sum_{j=1}^{m}\gamma_{j}$, the scalar $\eta\in(0,1)$ and $\delta>0$ satisfy $\delta \geq \frac{1}{m^{mB}},~~ \eta \leq(1-\frac{1}{m^{mB}})^{\frac{1}{mB}}.$
\end{theorem}

Theorem \ref{th3} provides that the bound of constraint violation measured by $||\sum_{i=1}^{m}A_{i}\widehat{\mathbf{x}}_{i}[T]-\mathbf{b}_{i}||$ is of
the order $O(\sqrt{\ln T/T})$.

\section{ Proof of main results} \label{4}

Before the proof of main results, we need to establish some useful auxiliary lemmas. The following Lemma \ref{le1} exploits the structure of strongly concave functions with Lipschitz gradients, whose proof is motivated by Lemma 3 in \cite{nedic2015distributed}. We omit the proof here.

\begin{lemma}\label{le1}
Let $h: \mathbb{R}^{p}\rightarrow \mathbb{R}$ be a $\gamma-$ strongly concave function
with $\gamma>0$ and have Lipschtiz continuous gradients with constant $M>0$. Let $\mathbf{z}\in \mathbb{R}^{p}$
and let $\mathbf{y}\in \mathbb{R}^{p}$ be defined by
$$\mathbf{y} = \mathbf{z} + \beta(\nabla h(\mathbf{z}) + \varphi(\mathbf{z})),$$
where $\beta\in (0, \frac{\gamma}{8M^{2}}]$ and $\varphi: \mathbb{R}^{p}\rightarrow \mathbb{R}^{p}$ is a mapping such that
$$||\varphi(\mathbf{z})||\leq c,~~~ \forall \mathbf{z}\in \mathbb{R}^{p}.$$
Then, there is a compact set $V \subset \mathbb{R}^{p}$ (which depends on $c$ and the define of function $h$, but not on $\beta$) such that
\begin{equation}
 ||\mathbf{y}|| \leq \left\{
\begin{aligned}
 ||\mathbf{z}||,~~ \forall \mathbf{z} \notin V,\\
  R   ,~~ \forall \mathbf{z} \in V ,\\
\end{aligned}
\right.\nonumber
\end{equation}
where $R = \max_{\mathbf{v}\in V}\{||\mathbf{v}|| + \frac{\gamma}{8M^{2}}||\nabla h(\mathbf{v})||\} + \frac{\gamma c}{8M^{2}}$.
\end{lemma}

Based on Lemma \ref{le1}, we are ready to prove our Theorem \ref{th1}.\\
{\bf Proof of Theorem \ref{th1}.} By step 5 of Algorithm \ref{alg1}, we have
$$\rho[t+1] = W[t]\rho[t],$$
where $\rho[t]$ is the vector with entries $\rho_{i}[t]$. Further, the above relation can be recursively written as follows
\begin{eqnarray}
\rho[t] = W[t-1]W[t-2]\cdots W[0]\mathbf{1},~~ \textrm{for}~~\textrm{all}~t \geq 1, \nonumber
\end{eqnarray}
where we use the fact that $\rho_{i}[0]=1$, for all $i = 1,2,\ldots, m$. Under Assumption \ref{A2}, by Corollary 2(b) in \cite{nedic2015distributed}, for all $i$, we have
\begin{eqnarray}
 \delta = \inf_{t = 0,1,\ldots} (\min_{1 \leq i \leq m}(W[t]W[t-1]\ldots W[0]\mathbf{1})_{i}) > 0. \nonumber
\end{eqnarray}
Therefore,  we can obtain
\begin{eqnarray}\label{eq-bound-p}
 \rho_{i}[t] \geq \delta, ~~\textrm{for}~~\textrm{all}~~i~~\textrm{and}~~t.
\end{eqnarray}
Using step 8 of Algorithm \ref{alg1}, we get
\begin{eqnarray}
\mathbf{\theta}_{i}[t] &= & \mathbf{u}_{i}[t] + \beta[t](A_{i}\mathbf{x}_{i}[t] - \mathbf{b}_{i} - \gamma_{i}\mathbf{\lambda}_{i}[t])\nonumber\\
&=& \rho_{i}[t] \left( \mathbf{\lambda}_{i}[t] + \frac{\beta[t]}{\rho_{i}[t]}(A_{i}\mathbf{x}_{i}[t] - \mathbf{b}_{i} - \gamma_{i}\mathbf{\lambda}_{i}[t]) \right),\nonumber
\end{eqnarray}
Furthermore, the above equality gives rise to
\begin{eqnarray}\label{eq-dual-variable-1}
\frac{\mathbf{\theta}_{i}[t]}{\rho_{i}[t]} &=& \mathbf{\lambda}_{i}[t] + \frac{\beta[t]}{\rho_{i}[t]}(A_{i}\mathbf{x}_{i}[t] - \mathbf{b}_{i}
- \gamma_{i}\mathbf{\lambda}_{i}[t])\nonumber\\
&=& (1 - \frac{\gamma_{i}\beta[t]}{\rho_{i}[t]})\mathbf{\lambda}_{i}[t] + \frac{\beta[t]}{\rho_{i}[t]}(A_{i}\mathbf{x}_{i}[t] - \mathbf{b}_{i}).
\end{eqnarray}

Since the transition matrix $W[t]W[t-1]\cdots W[0]$ is column stochastic and $\rho[0] = \mathbf{1}$, we have that $\sum_{i}^{m}\rho_{i}[t] = m$,
and $\delta \leq \rho_{i}[t] \leq m$. Together with (\ref{eq-bound-p}) and $\beta[t]\rightarrow 0$, it yields
\begin{eqnarray}
\lim_{t \rightarrow \infty} \frac{\mathbf{\beta}[t]}{\rho_{i}[t]}=0.
\end{eqnarray}
Thus, for each $i$, there exists a $T_{i} > 1$ such that $\frac{\beta[t]}{\rho_{i}[t]} \leq \frac{\tau_{i}^{3}}{8||A_{i}||^{2}}$, for all $t \geq T_{i}$.

Since the function $\phi_{i}(\lambda)$ defined in (\ref{dual}) is $\gamma_{i}$-strongly concave, and its gradient $\nabla\phi(\lambda)$ is Lipschitz continuous with constant $||A_{i}||/\tau_{i}$,
by Lemma \ref{le1}, there exists a finite $T_{i}>0$ and a compact set $V_{i}$ such that, for all $t \geq T_{i}$,
\begin{equation}\label{eq-dual-variable-2}
 \parallel\frac{\mathbf{\theta}_{i}[t]}{\rho_{i}[t]}\parallel \leq\left\{
\begin{aligned}
 ||\mathbf{\lambda}_{i}||,~~ \textrm{if}~~ \mathbf{\lambda}_{i} \notin V_{i},\\
  R_{i}   ,~~ \textrm{if}~~ \mathbf{\lambda}_{i} \in V_{i} .\\
\end{aligned}
\right.
\end{equation}

Let $T_{0}=\max_{1\leq i\leq m} T_{i}$. Now we divide $t$ into two part ($t\geq T_{0}$ and $1\leq t \leq T_{0}$) to prove the boundedness of $\parallel\frac{\mathbf{\theta}_{i}[t]}{\rho_{i}[t]}\parallel$, given by (\ref{eq-dual-variable-1}).

(i) By exploiting the mathematical induction, we will prove that, for all $t \geq T_{0}$,
\begin{eqnarray}\label{eq-bound-dual-variable}
\max_{1\leq i \leq m}||\lambda_{i}[t]|| \leq \widetilde{R},
\end{eqnarray}
where $\widetilde{R} = \max\{\max_{i}R_{i}, \max_{j}||\lambda_{j}[T]||\}$. Cleanly, if $t=T_{0}$, the relation (\ref{eq-bound-dual-variable}) is true.
Suppose it is true at some time $t \geq T_{0}$. Then, by (\ref{eq-dual-variable-2}), we have

\begin{eqnarray}\label{eq-dual-variable-3}
\parallel\frac{\mathbf{\theta}_{i}[t]}{\rho_{i}[t]}\parallel \leq  \max\left \{ R_{i}, \max_{j}||\mathbf{\lambda}_{j}[t]|| \right\} \leq
\widetilde{R},~~\textrm{for}~~ \textrm{all}~~i,
\end{eqnarray}
due to the induction hypothesis.

Next, in Lemma 4 of \cite{nedic2015stochastic}, we let $\mathbf{v} = \rho[t]$, $P = W[t]$, and $\mathbf{u}$ be taken as the vector of the $s$th coordinates of the vectors $\mathbf{\theta}_{i}[t]$, $i=1,\ldots,m$, where the coordinate index $s$ is arbitrary. By Lemma 4 of \cite{nedic2015stochastic}, we can get that each vector $\mathbf{\lambda}_{i}[t+1]$ is a convex combination of the vector $\frac{\mathbf{\theta}_{i}[t]}{\rho_{i}[t]}$, i.e.,
\begin{eqnarray}\label{eq-dual-variable-4}
\mathbf{\lambda}_{i}[t+1] = \sum_{j=1}^{m}Q_{ij}[t]\frac{\mathbf{\theta}_{i}[t]}{\rho_{i}[t]}, ~~\textrm{for}~~\textrm{all}~~i~~\textrm{and} ~~t\geq0,
\end{eqnarray}
where $Q[t]$ is a row stochastic matrix with entries $Q_{ij}[t] = \frac{W_{ij}[t]\rho_{j}[t]}{\rho_{i}[t+1]}$. Due to the convexity of Euclidean norm $||\cdot||$, we further obtain
\begin{eqnarray}\label{eq-bound-dual-variable-2}
||\mathbf{\lambda}_{i}[t+1]|| \leq \sum_{j=1}^{m}Q_{ij}[t]||\frac{\mathbf{\theta}_{j}[t]}{\rho_{j}[t]}|| \leq
 \max_{1\leq j\leq m} ||\frac{\mathbf{\theta}_{j}[t]}{\rho_{j}[t]}||, ~~\textrm{for}~~\textrm{all}~~i~~\textrm{and} ~~t\geq0.
\end{eqnarray}
By (\ref{eq-dual-variable-3}) and (\ref{eq-bound-dual-variable-2}), we have $||\mathbf{\lambda}_{i}[t+1]|| \leq \widetilde{R}$,
thus implying that, at time $t+1$
$$\max_{1\leq i\leq m}||\mathbf{\lambda}_{i}[t+1]|| \leq \widetilde{R}.$$
Hence, the relation (\ref{eq-dual-variable-3}) holds, for all $t\geq T_{0}$.

(ii) We prove that $||\mathbf{\lambda}_{i}[t]||$ is bounded upper when $t=1,2,\ldots, T_{0}-1$.
There is a constant $C>0$ such that $|1 - \frac{\gamma_{j}\beta[t]}{\rho_{j}[t]}|\leq C$, for all $t=1,2,\ldots, T_{0}-1$. Thus, together with (\ref{eq-dual-variable-1}) and (\ref{eq-dual-variable-4}), we can obtain that, for all $t=1,2,\ldots, T_{0}-1$
\begin{eqnarray}
\max_{1\leq i \leq m} ||\mathbf{\lambda}_{i}[t+1]|| &\leq& \max_{1\leq j\leq m} (1 - \frac{\gamma_{j}\beta[t]}{\rho_{j}[t]})||\mathbf{\lambda}_{j}[t]||
+ \frac{\beta[t]}{\rho_{j}[t]}||(A_{j}\mathbf{x}_{j}[t] - \mathbf{b}_{j})||\nonumber\\
&\leq &  \max_{1\leq j\leq m} C ||\mathbf{\lambda}_{j}[t]|| + \frac{\beta[t]}{\rho_{j}[t]}||(A_{j}\mathbf{x}_{j}[t] - \mathbf{b}_{j})||.\nonumber\\
&\leq &  \max_{1\leq j\leq m} C ||\mathbf{\lambda}_{j}[t]||+\max_{1\leq j\leq m}
\frac{\overline{\beta}}{\delta}G_{j},\nonumber
\end{eqnarray}
where $\overline{\beta}=\max_{1\leq t\leq T_{0}-1}\beta[t]$. Thus, exploiting the preceding relation recursively for $t=1,2,\ldots, T_{0}-1$,
and the fact that the initial point $\theta_{i}[0]$ is given in Algorithm \ref{alg1}, we conclude that there is a uniform deterministic bound on $||\mathbf{\lambda}_{i}[t]||$ for all $i$ and $t=1,2,\ldots, T_{0}-1$. According to the above discussion, we conclude the proof.  \qed

In order to prove Theorem \ref{th2} and \ref{th3}, we need to use the following result, which is a generalization of Lemma 8 in \cite{nedic2015distributed}.

\begin{lemma}\label{le2}
Under the conditions of Theorem \ref{th1}, for any $\mathbf{\lambda}\in \mathbb{R}^{d}$ and $t> 0$, we have
\begin{eqnarray}
||\overline{\mathbf{\theta}}[t+1] - \mathbf{\lambda}||^{2}
&\leq & ||\overline{\mathbf{\theta}}[t]-\mathbf{\lambda}||^{2} +
\frac{4\beta[t+1]}{m}\sum_{j=1}^{m}(G_{j} + \gamma_{j}D)||\mathbf{\lambda}_{j}[t+1] - \overline{\mathbf{\theta}}[t]||\nonumber\\
&& -\frac{\beta[t+1]}{m}\sum_{j=1}^{m}\gamma_{j}||\mathbf{\lambda}_{j}[t+1]-\mathbf{\lambda}||^{2} + \frac{\beta^{2}[t+1]}{m}\sum_{j=1}^{m}(G_{j}
 + \gamma_{j}D)^{2}\nonumber\\
&& -\frac{2\beta[t+1]}{m}(\mathcal{L}(\mathbf{x}[t+1],\mathbf{\lambda})-\mathcal{L}(\mathbf{x},\overline{\mathbf{\theta}}[t])).\nonumber
\end{eqnarray}
\end{lemma}
{\bf Proof }: We first prove that $\overline{\mathbf{\theta}}[t]$ is bounded, for any $t>0$. Since $W[t]$ is a column stochastic matrix, we have $\mathbf{1}^{\top}\mathbf{y}=\mathbf{1}^{\top}W[t]\mathbf{y}$, for
any vector $\mathbf{y}\in \mathbb{R}^{m}$.
By step 4 of Algorithm \ref{alg1}, we have
$$\sum_{i=1}^{m}\mathbf{u}_{i}[t+1]=\sum_{i=1}^{m}\sum_{j=1}^{m}W[t]_{ij}\mathbf{\theta}_{j}[t]=\sum_{j=1}^{m}\mathbf{\theta}_{j}[t] = m\overline{\theta}[t].$$
From the definition of $\mathbf{\lambda}_{i}[t+1]$ in step 6 of Algorithm \ref{alg1}, it gives rise to
\begin{eqnarray}\label{eq-proof-lemma2-1}
\overline{\theta}[t] = \frac{1}{m}\sum_{i=1}^{m}\mathbf{u}_{i}[t+1] = \frac{1}{m}\sum_{i=1}^{m}\rho_{i}[t+1]\mathbf{\lambda}_{i}[t+1].\nonumber
\end{eqnarray}
Note that $\sum_{i=1}^{m}\rho_{i}[t]=m$, and $\rho_{i}[t]>0$ for all $t$ and $i$. Thus, by the result of Theorem \ref{th1}, we have, for all $i=1,2,\ldots, m$ and $t\geq 0$,
\begin{eqnarray}\label{eq-proof-lemma2-2}
\overline{\theta}[t] \leq \max_{i}||\mathbf{\lambda}_{i}[t+1]|| \leq D.\nonumber
\end{eqnarray}

Now we are beginning to prove the result of Lemma \ref{le2}. From step 8 of Algorithm \ref{alg1}, we have
\begin{eqnarray}\label{eq-proof-lemma2-3}
\overline{\mathbf{\theta}}[t+1] = \overline{\mathbf{\theta}}[t] + \frac{\beta[t+1]}{m}\sum_{j=1}^{m}(A_{j}\mathbf{x}_{j}[t+1] -\mathbf{b}_{j} -
 \gamma_{j}\mathbf{\lambda}_{j}[t+1]).
\end{eqnarray}
For any $\mathbf{\lambda}\in \mathbb{R}^{d}$, the relation (\ref{eq-proof-lemma2-3}) gives rise to
\begin{eqnarray}
||\overline{\mathbf{\theta}}[t+1] -\mathbf{\lambda}||^{2}&=&||\overline{\mathbf{\theta}}[t]- \mathbf{\lambda} +
\frac{\beta[t+1]}{m}\sum_{j=1}^{m}(A_{j}\mathbf{x}_{j}[t+1] -\mathbf{b}_{j} - \gamma_{j}\mathbf{\lambda}_{j}[t+1])||^{2}\nonumber\\
\leq && ||\overline{\mathbf{\theta}}[t]- \mathbf{\lambda}||^{2}  + \frac{\beta^{2}[t+1]}{m^{2}}||\sum_{j=1}^{m}(A_{j}\mathbf{x}_{j}[t+1] -\mathbf{b}_{j} -
 \gamma_{j}\mathbf{\lambda}_{j}[t+1])||^{2}\nonumber\\
&&+ \frac{2\beta[t+1]}{m}\sum_{j=1}^{m}(A_{j}\mathbf{x}_{j}[t+1] -\mathbf{b}_{j} - \gamma_{j}\mathbf{\lambda}_{j}[t+1])^{\top}(\overline{\mathbf{\theta}}[t]-
 \mathbf{\lambda}).\nonumber
\end{eqnarray}
By using the inequality $(\sum_{j=1}^{m}a_{j})^{2} \leq m\sum_{j=1}^{m}a^{2}_{j}$, we can obtain
\begin{eqnarray}
||\sum_{j=1}^{m}(A_{j}\mathbf{x}_{j}[t+1] -\mathbf{b}_{j} - \gamma_{j}\mathbf{\lambda}_{j}[t+1]||^{2} &\leq & m\sum_{j=1}^{m}|| A_{j}\mathbf{x}_{j}[t+1]
-\mathbf{b}_{j} -
\gamma_{j}\mathbf{\lambda}_{j}[t+1]||^{2}\nonumber\\
&\leq & m\sum_{j=1}^{m}(G_{j}+\gamma_{j}D)^{2}.\nonumber
\end{eqnarray}
Thus, we have, for all $t\geq 0$
\begin{eqnarray}\label{eq-proof-lemma2-4}
&&||\overline{\mathbf{\theta}}[t+1] -\mathbf{\lambda}||^{2} \leq  ||\overline{\mathbf{\theta}}[t]- \mathbf{\lambda}||^{2} +
\frac{\beta^{2}[t+1]}{m}\sum_{j=1}^{m}(G_{j}+\gamma_{j}D)^{2}\nonumber\\
&&+ \frac{2\beta[t+1]}{m}\sum_{j=1}^{m}(A_{j}\mathbf{x}_{j}[t+1] -\mathbf{b}_{j} -
\gamma_{j}\mathbf{\lambda}_{j}[t+1])^{\top}(\overline{\mathbf{\theta}}[t]- \mathbf{\lambda}).
\end{eqnarray}
We now consider the last term in the right-hand side of (\ref{eq-proof-lemma2-4}), it can rewritten as
\begin{eqnarray}\label{eq-proof-lemma2-5}
&&(A_{j}\mathbf{x}_{j}[t+1] -\mathbf{b}_{j} - \gamma_{j}\mathbf{\lambda}_{j}[t+1])^{\top}(\overline{\mathbf{\theta}}[t]- \mathbf{\lambda})\nonumber\\
=&& (A_{j}\mathbf{x}_{j}[t+1] -\mathbf{b}_{j} - \gamma_{j}\mathbf{\lambda}_{j}[t+1])^{\top}(\overline{\mathbf{\theta}}[t]-\mathbf{\lambda}_{j}[t+1] +
\mathbf{\lambda}_{j}[t+1] - \mathbf{\lambda})\nonumber\\
=&& (A_{j}\mathbf{x}_{j}[t+1] -\mathbf{b}_{j} - \gamma_{j}\mathbf{\lambda}_{j}[t+1])^{\top}(\overline{\mathbf{\theta}}[t]-\mathbf{\lambda}_{j}[t+1])\nonumber\\
 &&+ (A_{j}\mathbf{x}_{j}[t+1] -\mathbf{b}_{j} - \gamma_{j}\mathbf{\lambda}_{j}[t+1])^{\top}(\mathbf{\lambda}_{j}[t+1] - \mathbf{\lambda}).
\end{eqnarray}
By the Cauchy-Schwarz inequality, we have
\begin{eqnarray}\label{eq-proof-lemma2-6}
&&-(G_{j}+\gamma_{j}D)||\overline{\mathbf{\theta}}[t]-\mathbf{\lambda}_{j}[t+1]|| \nonumber\\
&&\leq (A_{j}\mathbf{x}_{j}[t+1] -\mathbf{b}_{j} - \gamma_{j}\mathbf{\lambda}_{j}[t+1])^{\top}(\overline{\mathbf{\theta}}[t]-\mathbf{\lambda}_{j}[t+1]).
\end{eqnarray}
Since $\mathcal{L}_{j}(\mathbf{x},\cdot)$ is $\gamma_{j}$-strongly concave, we have, for any $\lambda \in \mathbb{R}^{p}$
\begin{eqnarray}\label{eq-proof-lemma2-7}
&&(A_{j}\mathbf{x}_{j}[t+1] -\mathbf{b}_{j} - \gamma_{j}\mathbf{\lambda}_{j}[t+1])^{\top}(\mathbf{\lambda}_{j}[t+1] - \mathbf{\lambda})\nonumber\\
\leq && \mathcal{L}_{j}(\mathbf{x}_{j}[t+1],\mathbf{\lambda}_{j}[t+1]) - \mathcal{L}_{j}(\mathbf{x}_{j}[t+1],\mathbf{\lambda}) -
\frac{\gamma_{j}}{2}||\mathbf{\lambda}_{j}[t+1] - \mathbf{\lambda}||^{2}.
\end{eqnarray}
By step 7 of Algorithm \ref{alg1}, for any $\mathbf{x}_{j} \in \mathbb{R}^{n}$, we can get
\begin{eqnarray}\label{eq-proof-lemma2-8}
&&f_{j}(\mathbf{x}_{j}[t+1]) + \mathbf{\lambda}_{j}[t+1]^{\top}(A_{j}\mathbf{x}_{j}[t+1] -\mathbf{b}_{j}) -
 \frac{\gamma_{j}}{2}\mathbf{\lambda}_{j}[t+1]^{\top}\mathbf{\lambda}_{j}[t+1]\nonumber\\
\leq&&  f_{j}(\mathbf{x}_{j}) + \mathbf{\lambda}_{j}[t+1]^{\top}(A_{j}\mathbf{x}_{j} -\mathbf{b}_{j}) - \frac{\gamma_{j}}{2}\mathbf{\lambda}_{j}[t+1]^{\top}\mathbf{\lambda}_{j}[t+1]. \nonumber
\end{eqnarray}
Subtracting $\mathcal{L}_{j}(\mathbf{x}_{j}[t+1],\mathbf{\lambda})$ in above relation, we obtain
\begin{eqnarray}\label{eq-proof-lemma2-9}
&&\mathcal{L}_{j}(\mathbf{x}_{j}[t+1],\mathbf{\lambda}_{j}[t+1]) - \mathcal{L}_{j}(\mathbf{x}_{j}[t+1],\mathbf{\lambda})\nonumber\\
&\leq& f_{j}(\mathbf{x}_{j}) + \mathbf{\lambda}_{j}[t+1]^{\top}(A_{j}\mathbf{x}_{j} -\mathbf{b}_{j}) -
 \frac{\gamma_{j}}{2}\mathbf{\lambda}_{j}[t+1]^{\top}\mathbf{\lambda}_{j}[t+1] - \mathcal{L}_{j}(\mathbf{x}_{j}[t+1],\mathbf{\lambda})\nonumber\\
&\leq& (G_{j}+\gamma_{j}D)||\mathbf{\lambda}_{j}[t+1]-\overline{\mathbf{\theta}}[t]|| +
 \mathcal{L}_{j}(\mathbf{x}_{j},\overline{\mathbf{\theta}}[t]) - \mathcal{L}_{j}(\mathbf{x}_{j}[t+1],\mathbf{\lambda}).
\end{eqnarray}
Together with (\ref{eq-proof-lemma2-5}), (\ref{eq-proof-lemma2-6}), (\ref{eq-proof-lemma2-7}), (\ref{eq-proof-lemma2-9}) and the definition of $\mathcal{L}(\mathbf{x},\lambda)$, we can obtain the desired result.  \qed

Next, we prove Theorem \ref{th2}.\\
{\bf Proof of Theorem \ref{th2}}~~Let $\mathbf{x}=\mathbf{x}^{*}$ and $\mathbf{\lambda}=\mathbf{0}$ in Lemma \ref{le2}, we have
\begin{eqnarray}\label{eq-proof-theorem3-1}
&&||\overline{\mathbf{\theta}}[t+1]-\mathbf{0}||^{2} \nonumber\\
\leq &&||\overline{\mathbf{\theta}}[t]-\mathbf{0}||^{2} + \frac{4\beta[t+1]}{m}\sum_{j=1}^{m}(G_{j}+\gamma_{j}D)||\lambda_{j}[t+1] -\overline{\theta}[t]||\nonumber\\
&&- \frac{2\beta[t+1]}{m}(\mathcal{L}(\mathbf{x}[t+1],\mathbf{0}) - \mathcal{L}(\mathbf{x}^{*},\overline{\mathbf{\theta}}[t])) -
\frac{\beta[t+1]}{m}\sum_{j=1}^{m}\gamma_{j}||\lambda_{j}[t+1]-\mathbf{0}||^{2}\nonumber\\
&&+ \frac{\beta^{2}[t+1]}{m}\sum_{j=1}^{m}(G_{j}+\gamma_{j}D)^{2}\nonumber\\
\leq && ||\overline{\mathbf{\theta}}[t]-\mathbf{0}||^{2}+ \frac{4\beta[t+1]}{m}\sum_{j=1}^{m}(G_{j}+\gamma_{j}D)||\lambda_{j}[t+1] -\overline{\theta}[t]||\nonumber\\
&&- \frac{2\beta[t+1]}{m}(\mathcal{L}(\mathbf{x}[t+1],\mathbf{0}) - \mathcal{L}(\mathbf{x}^{*},\overline{\mathbf{\theta}}[t]))+
 \frac{\beta^{2}[t+1]}{m}\sum_{j=1}^{m}(G_{j}+\gamma_{j}D)^{2}.
\end{eqnarray}
Using the definition of function $\mathcal{L}(\mathbf{x},\lambda)$ and letting $ \gamma= \sum_{j=1}^{m}\gamma_{j}$, we can obtain
\begin{eqnarray}\label{eq-proof-theorem3-2}
&&\mathcal{L}(\mathbf{x}[t+1],\mathbf{0}) - \mathcal{L}(\mathbf{x}^{*},\overline{\mathbf{\theta}}[t])\nonumber\\
=&& \mathcal{L}(\mathbf{x}[t+1],\mathbf{0}) - \mathcal{L}(\mathbf{x}^{*},\mathbf{0}) + \mathcal{L}(\mathbf{x}^{*},\mathbf{0})-
 \mathcal{L}(\mathbf{x}^{*},\overline{\mathbf{\theta}}[t])\nonumber\\
=&& F(\mathbf{x}[t+1])-F(\mathbf{x}^{*}) + \mathcal{L}(\mathbf{x}^{*},\mathbf{0})- \mathcal{L}(\mathbf{x}^{*},\overline{\mathbf{\theta}}[t])\nonumber\\
\geq && F(\mathbf{x}[t+1])-F(\mathbf{x}^{*}) + \frac{\gamma}{2}||\mathbf{\overline{\theta}}[t]-\mathbf{0}||^{2},
\end{eqnarray}
where the last inequality makes use of the strong concavity of $\mathcal{L}(\mathbf{x},\cdot)$. Thus, by (\ref{eq-proof-theorem3-1}) and (\ref{eq-proof-theorem3-2}), and then letting $\beta[t]=\frac{q}{t}$, we have
\begin{eqnarray}
&&||\overline{\mathbf{\theta}}[t+1]-\mathbf{0}||^{2} \nonumber\\
\leq && (1-\frac{q\gamma}{m(t+1)})||\overline{\mathbf{\theta}}[t]-\mathbf{0}||^{2}+
\frac{4q}{m(t+1)}\sum_{j=1}^{m}(G_{j}+\gamma_{j}D)||\lambda_{j}[t+1] -\overline{\theta}[t]||\nonumber\\
&&-\frac{2q}{m(t+1)}(F(\mathbf{x}[t+1])-F(\mathbf{x}^{*})) +\frac{q^{2}}{m(t+1)^{2}}\sum_{j=1}^{m}(G_{j}+\gamma_{j}D)^{2}.\nonumber
\end{eqnarray}
Note that $4\leq\frac{q\gamma}{m}$, it follows that
\begin{eqnarray}
&&||\overline{\mathbf{\theta}}[t+1]-\mathbf{0}||^{2} \nonumber\\
\leq && (1-\frac{2}{t+1})||\overline{\mathbf{\theta}}[t]-\mathbf{0}||^{2}+ \frac{4q}{m(t+1)}\sum_{j=1}^{m}(G_{j}+\gamma_{j}D)||\lambda_{j}[t+1] -\overline{\theta}[t]||
 \nonumber\\
&&-\frac{2q}{m(t+1)}(F(\mathbf{x}[t+1])-F(\mathbf{x}^{*})) +\frac{q^{2}}{m(t+1)^{2}}\sum_{j=1}^{m}(G_{j}+\gamma_{j}D)^{2}.\nonumber
\end{eqnarray}
Multiplying the preceding relation by $t(t+1)$, we can see that, for all $t\geq1$
\begin{eqnarray}\label{eq-proof-theorem3-3}
&&(t+1)t||\overline{\mathbf{\theta}}[t+1]-\mathbf{0}||^{2} \nonumber\\
\leq && t(t-1)||\overline{\mathbf{\theta}}[t]-\mathbf{0}||^{2}+ \frac{4qt}{m}\sum_{j=1}^{m}(G_{j}+\gamma_{j}D)||\lambda_{j}[t+1] -\overline{\theta}[t]|| \nonumber\\
&&-\frac{2qt}{m}(F(\mathbf{x}[t+1])-F(\mathbf{x}^{*})) +\frac{q^{2}t}{m(t+1)}\sum_{j=1}^{m}(G_{j}+\gamma_{j}D)^{2}.\nonumber
\end{eqnarray}
Summing up the above inequality from $1$ to $(T-1)$ for all $T\geq2$ and rearranging the terms, it leads to
\begin{eqnarray}\label{eq-proof-theorem3-4}
&&\frac{2q}{m}\sum_{t=1}^{T-1}t(F(\mathbf{x}[t+1])-F(\mathbf{x}^{*})) \nonumber\\
\leq && -T(T-1)||\overline{\mathbf{\theta}}[t]-\mathbf{0}||^{2} + \frac{q^{2}}{m}\sum_{t=1}^{T-1}\frac{t}{t+1}\sum_{j=1}^{m}(G_{j}+\gamma_{j}D)^{2}\nonumber\\
&& + \frac{4q}{m}\sum_{t=1}^{T-1}t\sum_{j=1}^{m}(G_{j}+\gamma_{j}D)||\lambda_{j}[t+1] -\overline{\theta}[t]||\nonumber\\
\leq && \frac{q^{2}(T-1)}{m}\sum_{j=1}^{m}(G_{j}+\gamma_{j}D)^{2} + \frac{4q}{m}\sum_{t=1}^{T-1}t\sum_{j=1}^{m}(G_{j}+
\gamma_{j}D)||\lambda_{j}[t+1] -\overline{\theta}[t]||\nonumber\\
\leq && \frac{q^{2}(T-1)}{m}\sum_{j=1}^{m}(G_{j}+\gamma_{j}D)^{2} + \frac{4q(T-1)}{m}\sum_{t=1}^{T-1}\sum_{j=1}^{m}(G_{j}+
\gamma_{j}D)||\lambda_{j}[t+1] -\overline{\theta}[t]||.\nonumber
\end{eqnarray}
Dividing both sides by $\frac{qT(T-1)}{m}$ in above relation, it yields
\begin{eqnarray}\label{eq-proof-theorem3-5}
&&\frac{2}{T(T-1)}\sum_{t=1}^{T-1}t(F(\mathbf{x}[t+1])-F(\mathbf{x}^{*})) \nonumber\\
\leq && \frac{4}{T}\sum_{t=1}^{T-1}\sum_{j=1}^{m}(G_{j}+\gamma_{j}D)||\lambda_{j}[t+1] -\overline{\theta}[t]||+\frac{q}{T}\sum_{j=1}^{m}(G_{j}+\gamma_{j}D)^{2}.
\end{eqnarray}

Note that, for all $i$ and $t$, we get
$$||A_{i}\mathbf{x}_{i}[t+1]-\mathbf{b}_{i}-\gamma_{i}\mathbf{\lambda}_{i}[t+1]||_{1}\leq \sqrt{p}||(A_{i}\mathbf{x}_{i}[t+1]-\mathbf{b}_{i}-\
\gamma_{i}\mathbf{\lambda}_{i}[t+1])||\leq \sqrt{p}(G_{i}+\gamma_{i}D).$$
Letting $\mathbf{e}_{i}[t]=\beta[t](A_{i}\mathbf{x}_{i}[t+1]-\mathbf{b}_{i}-\gamma_{i}\mathbf{\lambda}_{i}[t+1])$ with $\beta[t]=\frac{q}{t}$, we have $||\mathbf{e}_{i}[t]||_{1}\leq \frac{q B}{t}$ for all $i$ and $t$, where $B=\max_{1\leq i\leq m}\sqrt{p}(G_{i}+\gamma_{i}D)$.
By applying Corollary 2 in \cite{nedic2015stochastic}, we can estimate the term $||\lambda_{j}[t+1] - \overline{\mathbf{\theta}}[t]||$ in (\ref{eq-proof-theorem3-5}) as follows
\begin{eqnarray}\label{eq-proof-theorem2-9}
\sum_{t=1}^{T-1}||\lambda_{j}[t+1] - \overline{\mathbf{\theta}}[t]||\leq\frac{8\eta}{\delta(1-\eta)}\sum_{j=1}^{m}||\mathbf{\theta}_{j}[0]||_{1} + \frac{8 q m B}{\delta(1-\eta)}(1+\ln T).
\end{eqnarray}
Combining  (\ref{eq-proof-theorem2-9}) with (\ref{eq-proof-theorem3-5}), we can get
\begin{eqnarray}\label{eq-proof-theorem3-6}
&&\frac{\sum_{t=1}^{T-1}t \left (F(\mathbf{x}[t+1])-F(\mathbf{x}^{*})\right )}{\frac{T(T-1)}{2}}\nonumber\\
\leq&&\frac{4}{T\delta}\sum_{j=1}^{m}(G_{j}+\gamma_{j}D)\left (\frac{8\eta}{1-\eta}\sum_{j=1}^{m}||\mathbf{\theta}_{j}[0]||_{1}+
\frac{8q m B}{1-\eta}(1+\ln T)\right ) \nonumber\\
&& + \frac{q}{T}\sum_{j=1}^{m}(G_{j}+\gamma_{j}D)^{2}.
\end{eqnarray}
Using the convexity of $F$, the definition of $\widehat{\mathbf{x}}[T]$ and (\ref{eq-proof-theorem3-6}), the desired result can be obtained.   \qed
\\
{\bf Proof of Theorem 3} Let $\mathbf{x}=\mathbf{x}^{*}$ in Lemma \ref{le2}, we can get
\begin{eqnarray}\label{eq-proof-theorem4-1}
&&\frac{2\beta[t+1]}{m}(\mathcal{L}(\mathbf{x}[t+1],\lambda) - \mathcal{L}(\mathbf{x}^{*},\overline{\mu}[t]))\nonumber\\
&&\leq ||\overline{\mu}[t]-\lambda||^{2}-||\overline{\mu}[t+1]-\lambda||^{2} + \frac{4\beta[t+1]}{m}\sum_{j=1}^{m}(G_{j} + \gamma_{j}D)||\lambda_{j}[t+1] - \overline{\mu}[t]||\nonumber\\
&&+ \frac{\beta^{2}[t+1]}{m}\sum_{j=1}{m}(G_{j}+\gamma_{j}D)^{2}.
\end{eqnarray}
Considering the terms in the left-hand side of (\ref{eq-proof-theorem4-1}), we have
\begin{eqnarray}\label{eq-proof-theorem4-2}
&&2(\mathcal{L}(\mathbf{x}[t+1],\lambda)-\mathcal{L}(\mathbf{x}^{*},\overline{\mu}[t]))\nonumber\\
&&=F(\mathbf{x}[t+1]) + \lambda^{\top}(\sum_{j=1}^{m}A_{j}\mathbf{x}_{j}[t+1]-\mathbf{b}_{j})-F(\mathbf{x}^{*}) - \frac{\gamma}{2}\lambda^{\top}\lambda+\frac{\gamma}{2}\overline{\mu}^{\top}[t]\overline{\mu}[t]\nonumber\\
&& +\mathcal{L}(\mathbf{x}[t+1],\lambda) - \mathcal{L}(\mathbf{x}^{*},\overline{\mu}[t])\nonumber\\
&&=F(\mathbf{x}[t+1]) + \lambda^{\top}(\sum_{j=1}^{m}A_{j}\mathbf{x}_{j}[t+1]-\mathbf{b}_{j})-F(\mathbf{x}^{*}) - \frac{\gamma}{2}\lambda^{\top}\lambda + \frac{\gamma}{2}\overline{\mu}^{\top}[t]\overline{\mu}[t]\nonumber\\
&&+ \mathcal{L}(\mathbf{x}[t+1],\lambda)-\mathcal{L}(\mathbf{x}^{*},\lambda)+\mathcal{L}(\mathbf{x}^{*},\lambda) - \mathcal{L}(\mathbf{x}^{*},\overline{\mu}[t])\nonumber\\
&&\geq 2(F(\mathbf{x}[t+1]) + \lambda^{\top}(\sum_{j=1}^{m}A_{j}\mathbf{x}_{j}[t+1]-\mathbf{b}_{j})-F(\mathbf{x}^{*}))- \frac{\gamma}{2}\lambda^{\top}\lambda + \frac{\gamma}{2}\overline{\mu}^{\top}[t]\overline{\mu}[t]\nonumber\\
&& + \frac{\gamma}{2}||\overline{\mu}[t]-\lambda||^{2}+\gamma\lambda^{\top}(\overline{\mu}[t]-\lambda),
\end{eqnarray}
where  the last inequality is due to the strong concavity of $L(\mathbf{x}^{*},\cdot)$. Further, by (\ref{eq-proof-theorem4-2}), we can deduce
\begin{eqnarray}\label{eq-proof-theorem4-4}
&& 2(\mathcal{L}(\mathbf{x}[t+1],\lambda)-\mathcal{L}(\mathbf{x}^{*},\overline{\mu}[t]))\nonumber\\
&&\geq 2\lambda^{\top}(\sum_{j=1}^{m}A_{j}\mathbf{x}_{j}[t+1]-\mathbf{b}_{j}) + \frac{\gamma}{2}||\overline{\mu}[t]-\lambda||^{2} -\frac{\gamma}{2}(3\lambda^{\top}\lambda - \overline{\mu}^{\top}[t]\overline{\mu}[t]-2\lambda^{\top}\overline{\mu}[t])\nonumber\\
&&= 2\lambda^{\top}(\sum_{j=1}^{m}A_{j}\mathbf{x}_{j}[t+1]-\mathbf{b}_{j}) + \frac{\gamma}{2}||\overline{\mu}[t]-\lambda||^{2} -\frac{\gamma}{2}(4\lambda^{\top}\lambda - ||\lambda+\overline{\mu}[t]||^{2}) \nonumber\\
&&\geq  2\lambda^{\top}(\sum_{j=1}^{m}A_{j}\mathbf{x}_{j}[t+1]-\mathbf{b}_{j}) + \frac{\gamma}{2}||\overline{\mu}[t]-\lambda||^{2}+2\gamma\lambda^{\top}\lambda.
\end{eqnarray}
Combining (\ref{eq-proof-theorem4-1}) with (\ref{eq-proof-theorem4-4}), and then letting $\beta[t+1]=\frac{q}{t+1}$, we can obtain
\begin{eqnarray}\label{eq-proof-theorem4-5}
&& \frac{2q}{m(t+1)}\left(\lambda^{\top}(\sum_{j=1}^{m}A_{j}\mathbf{x}_{j}[t+1]-\mathbf{b}_{j})-\gamma\lambda^{\top}\lambda\right)\nonumber\\
&&\leq (1-\frac{q\gamma}{2m(t+1)})||\overline{\mu}[t]-\lambda||^{2} - ||\overline{\mu}[t+1]-\lambda||^{2} + \frac{4q}{m(t+1)}\sum_{j=1}^{m}(G_{j}+ \gamma_{j}D)||\lambda[t+1]-\overline{\mu}[t]||\nonumber\\
&& + \frac{q^{2}}{m(t+1)^{2}}\sum_{j=1}^{2}(G_{j}+ \gamma_{j}D)^{2}.\nonumber
\end{eqnarray}
Due to the fact that $4 \leq \frac{q\gamma}{m}$, we can see that $1-\frac{q\gamma}{2m(t+1)} \leq 1-\frac{2}{t+1}$. Thus, by the preceding inequality,  we can obtain
\begin{eqnarray}\label{eq-proof-theorem4-6}
&& \frac{2q}{m(t+1)}\left(\lambda^{\top}(\sum_{j=1}^{m}A_{j}\mathbf{x}_{j}[t+1]-\mathbf{b}_{j})-\gamma\lambda^{\top}\lambda\right)\nonumber\\
&&\leq (1-\frac{2}{t+1})||\overline{\mu}[t]-\lambda||^{2} - ||\overline{\mu}[t+1]-\lambda||^{2} + \frac{4q}{m(t+1)}\sum_{j=1}^{m}(G_{j}+ \gamma_{j}D)||\lambda[t+1]-\overline{\mu}[t]||\nonumber\\
&& + \frac{q^{2}}{m(t+1)^{2}}\sum_{j=1}^{2}(G_{j}+ \gamma_{j}D)^{2}.\nonumber
\end{eqnarray}
Multiplying the above inequality by $t(t+1)$, and then summing up from $1$ to $T-1$, we have, for all $t\geq 1$ and $ T\geq 2$

\begin{eqnarray}\label{eq-proof-theorem4-8}
&& \frac{2q}{m}\sum_{t=1}^{T-1}t\left(\lambda^{\top}(\sum_{j=1}^{m}A_{j}\mathbf{x}_{j}[t+1]-\mathbf{b}_{j})-\gamma\lambda^{\top}\lambda\right)\nonumber\\
&&\leq  \frac{4q}{m}\sum_{t=1}^{T-1}t\sum_{j=1}^{m}(G_{j}+ \gamma_{j}D)||\lambda[t+1]-\overline{\mu}[t]|| + \frac{q^{2}}{m}\sum_{t=1}^{T-1}\frac{t}{t+1}\sum_{j=1}^{2}(G_{j}+ \gamma_{j}D)^{2}\nonumber\\
&& -T(T-1)||\overline{\mu}[T]-\lambda||^{2}\nonumber\\
&& \leq \frac{4q(T-1)}{m}\sum_{t=1}^{T-1}\sum_{j=1}^{m}(G_{j}+ \gamma_{j}D)||\lambda[t+1]-\overline{\mu}[t]|| + \frac{q^{2}(T-1)}{m}\sum_{j=1}^{2}(G_{j}+ \gamma_{j}D)^{2}.\nonumber
\end{eqnarray}
Dividing both sides by $\frac{qT(T-1)}{m}$ in the inequality above, it gives rise to
\begin{eqnarray}\label{eq-proof-theorem4-9}
&& \frac{\sum_{t=1}^{T-1}t(\lambda^{\top}(\sum_{j=1}^{m}A_{j}\mathbf{x}_{j}[t+1]-\mathbf{b}_{j})-\gamma\lambda^{\top}\lambda)}{\frac{(T-1)T}{2}}\nonumber\\
&& \leq \frac{q}{T}\sum_{j=1}^{m}(G_{j}+ \gamma_{j}D)^{2} + \frac{4}{T}\sum_{t=1}^{T-1}\sum_{j=1}^{m}(G_{j}+ \gamma_{j}D)||\lambda[t+1]-\overline{\mu}[t]||.
\end{eqnarray}
Note that $\sum_{j=1}^{m}A_{j}\mathbf{x}_{j}[t+1]-\mathbf{b}_{j}$ is linear, thus, we have
\begin{eqnarray}\label{eq-proof-theorem4-10}
&& \frac{\sum_{t=1}^{T-1}t(\lambda^{\top}(\sum_{j=1}^{m}A_{j}\mathbf{x}_{j}[t+1]-\mathbf{b}_{j})-\gamma\lambda^{\top}\lambda)}{\frac{(T-1)T}{2}} \geq \lambda^{\top}(\sum_{j=1}^{m}A_{j}\widehat{\mathbf{x}}_{j}[T]-\mathbf{b}_{j})-\gamma\lambda^{\top}\lambda.
\end{eqnarray}
By (\ref{eq-proof-theorem4-9}) and (\ref{eq-proof-theorem4-10}), we can obtain, for any $\lambda\in \mathbb{R}^{p}$
\begin{eqnarray}\label{eq-proof-theorem4-11}
&& \lambda^{\top}(\sum_{j=1}^{m}A_{j}\widehat{\mathbf{x}}_{j}[T]-\mathbf{b}_{j})-\gamma\lambda^{\top}\lambda\nonumber\\
&& \leq \frac{q}{T}\sum_{j=1}^{m}(G_{j}+ \gamma_{j}D)^{2} + \frac{4}{T}\sum_{t=1}^{T-1}\sum_{j=1}^{m}(G_{j}+ \gamma_{j}D)||\lambda[t+1]-\overline{\mu}[t]||.
\end{eqnarray}
Maximizing the terms in the left-hand side of (\ref{eq-proof-theorem4-11}) with respect to $\lambda$  and using the estimate (\ref{eq-proof-theorem2-9}), we can get the desired result.
The proof is completed.            \qed

\section{Numerical experiments}\label{5}

Distributed optimization problems with coupled equality constraints have an interesting application on the network utility maximization (NUM) problem investigated in \cite{Low1999optimization,beck2014an,necoara2008application}.
More specifically, a network is modeled as a set of links $L$ with finite capacities $C=(C_l,l\in L)$. They are shared by a set of sources $S$
indexed by $s$. Each source $s$ uses a set $L(s)\subset L$. Let $S(l)=\{s\in S| l\in L(s)\}$ be the set of sources using
link $l$. The set $\{L(s)\}$ defines an $|L|\times |S|$ routing matrix $A$ with entries given by $A_{ls}=1$ if $l\in L(s)$, $A_{ls}=0$ otherwise.
Each source $s$ is associated with a utility function $U_s: \mathbb{R}^+ \rightarrow  \mathbb{R}$, i.e., source $s$ gains a utility $U_s(x_s)$ when it sends data at rate $x_s$ satisfying
$0\leq m_s\leq x_s\leq M_s$. Let $I_s=[m_s,M_s]$.
Mathematically, the NUM problem is to determine the source rates that minimize the sum of disutilities with link capacity constraints \cite{Low1999optimization}:
\begin{eqnarray}
\textrm{(NUM)}~~~\min\limits_{x_s\in I_s} && g_N(x):=\sum_{s\in S}-U_{s}(x_s)  \nonumber\\
  \textrm{s.t.} && Ax= C.   \nonumber
\end{eqnarray}

Note that the utility function $U_s$ and constraint $I_s$ are \emph{local} and \emph{private}, only known by the source $s$. Solving the NUM problem directly requires coordination among possibly all sources and is impractical in real networks. It is important to seek a distributed solution. In the following numerical experiments, we will utilize our proposed distributed method to solve the NUM problem.

For numerical simulations, the utility function is taken as $U_s(x_s)=20 w_s \log(x_s+0.1)$ from \cite{beck2014an}.
Set $C_l=1$ for all $l\in L$, and $w_s={|L(s)|}/{|L|},m_s=0,M_s=1$ for all $s\in S$. For the communicated weight matrix $W[t]$, a pool of 20 weight matrices connecting random graphs are generated, in which each weight matrix satisfies Assumption \ref{A2}. We take all the regularization parameters as the same with $\gamma_{s}=1, s\in S$ and the stepsize parameter as $q=4$. We use MATLAB convex programming
toolbox CVX  to compute the solution $x^{*}$. For our method and the compared algorithm, all the algorithms were terminated when all of the conditions below are satisfied at an iteration $t$:
(i) $\max_{s\in S}|\lambda_s[t+1]-\lambda_s[t]|\leq \epsilon$,
(ii) $\max_{l\in L} ||Ax[t+1]-C|| \leq \epsilon$,
(iii) $\max_{s\in S}|\frac{U_{s}(x[t+1])-U_{s}(x[t])}{U_{s}(x[t])}| \leq \epsilon$,
where we set $\epsilon=0.01$ in the simulations.

\begin{figure}[!htbp]
\centering
\includegraphics[width=10cm]{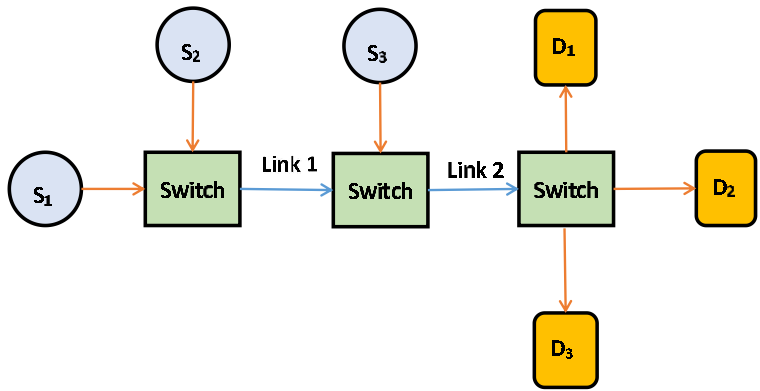}
\caption{Logical topology. Source $S_{i},i=1,2,3$ transmits to destinations $D_{i},i=1,2,3$}
\label{fig7}
\end{figure}

\begin{figure}[!htbp]
\centering
\includegraphics[width=10cm]{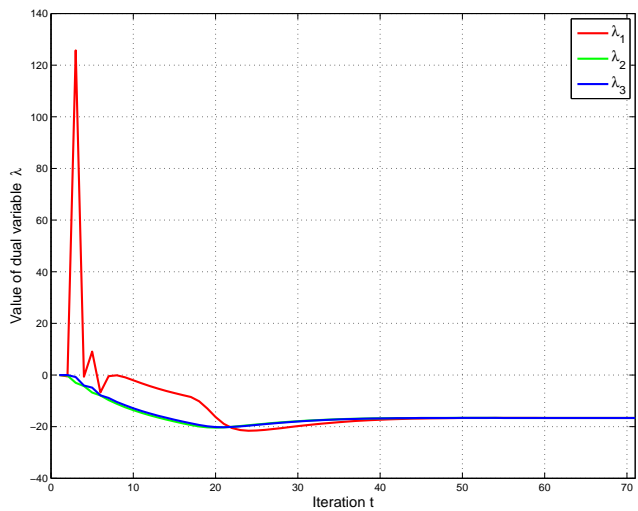}
\caption{Iterative value of dual variable $\lambda$}
\label{fig1}
\end{figure}

\begin{figure}[!htbp]
\centering
\includegraphics[width=12cm]{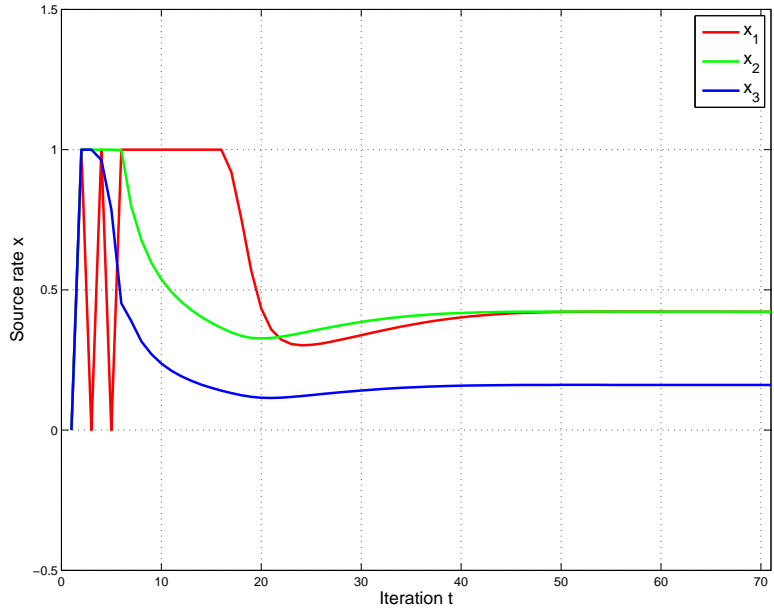}
\caption{Iterative value of source rate $x$}
\label{fig2}
\end{figure}

\begin{figure}[!htbp]
\centering
\includegraphics[width=12cm]{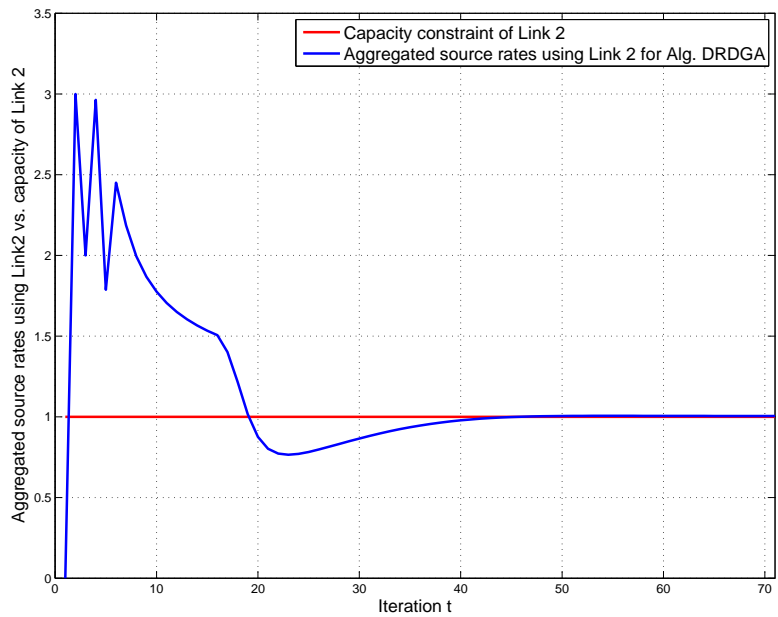}
\caption{Aggregated source rates using Link 2 vs. capacity of Link 2}
\label{fig3}
\end{figure}

\begin{figure}[!htbp]
\centering
\includegraphics[width=12cm]{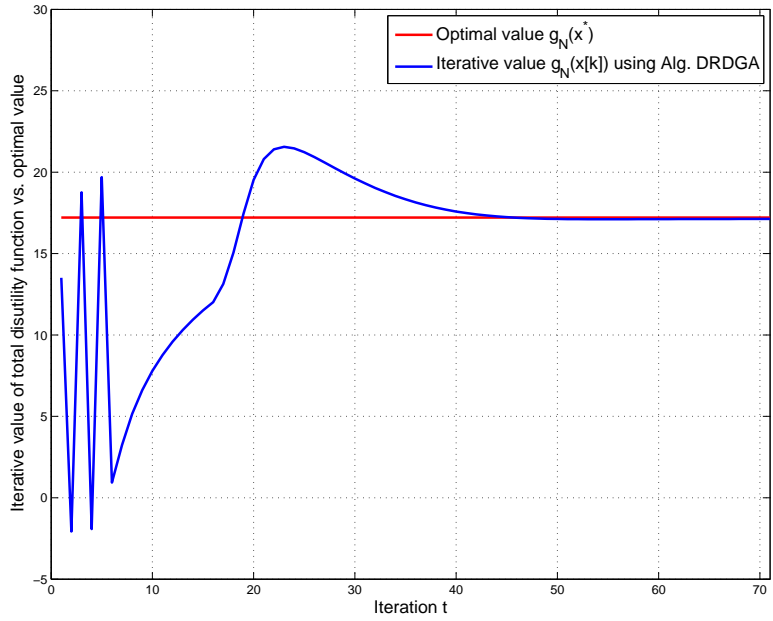}
\caption{Iterative value of total disutity function vs. optimal value}
\label{fig4}
\end{figure}

We first consider a simple logical topology with $S=3$ and $L=2$ \cite{Low1999optimization}, displayed as in  Figure \ref{fig7}. It follows from Figure \ref{fig7} that $w_{1}=1, w_{2}=1, w_{3}=1/2$. Figure \ref{fig1} shows the evolution of dual variables at the first 70 iterations. Clearly, all local dual variables $\lambda_{s}, s=1,2,3$, agree on the same value at a short time with around 70 iterations. Figure \ref{fig2} illustrates the evolution of each source rate $x_{s}, s=1,2,3$. Source rate $x_{1}$ and $x_{2}$ can arrive at same value because the weight coefficients $w_{1}=w_{2}$. After 70 iterations, every source rate $x_{s}$ can arrive approximately at the optimal solution. Figure \ref{fig3} demonstrates the aggregated source rates that use Link 2 versus capacity limit of Link 2. It can observed from Figure \ref{fig3} that the aggregated source rates satisfy the constraint of Link 2 capacity appropriately. As shown in Figure \ref{fig4}, the iterative values of disutility objective function $g_{N}(x[t])$ rapidly converge to the optimal value $g_{N}(x^{*})$ .

\begin{figure}[!htbp]
\centering
\includegraphics[width=12cm]{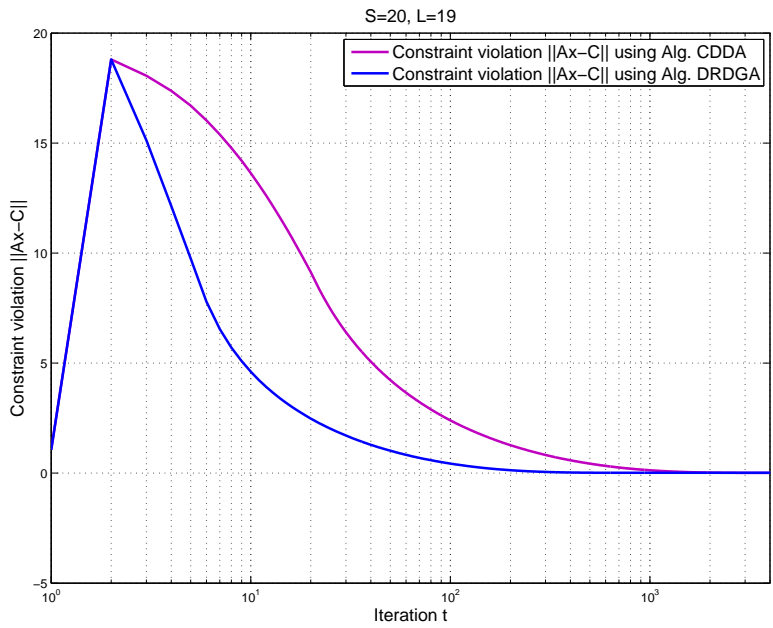}
\caption{Evolution of constraint violation $||Ax[t]-C||$}
\label{fig5}
\end{figure}

To compare the performance of our proposed Alg. DRDGA with the existing dual decomposition distributed algorithm (Alg. CDDA) in \cite{falsone2016dual}, we next test a random generated problem NUM with sizes $S=20, L=19$ and report the comparisons on the constraint violations and objective function values. Figure \ref{fig5} displays the evolution of the constraint violation $||Ax[t]-C||$. We can find that both algorithms can satisfy the linear equality constraints gradually. But, the convergence speedup of our Alg. DRDGA is faster than that of Alg. CDDA. Figure \ref{fig5} illustrates that both algorithms can also converge to the optimal value. However, by comparisons, our Alg. DRDGA is convergent to the optimal value faster than Alg. CDDA.

\begin{figure}[!htbp]
\centering
\includegraphics[width=12cm]{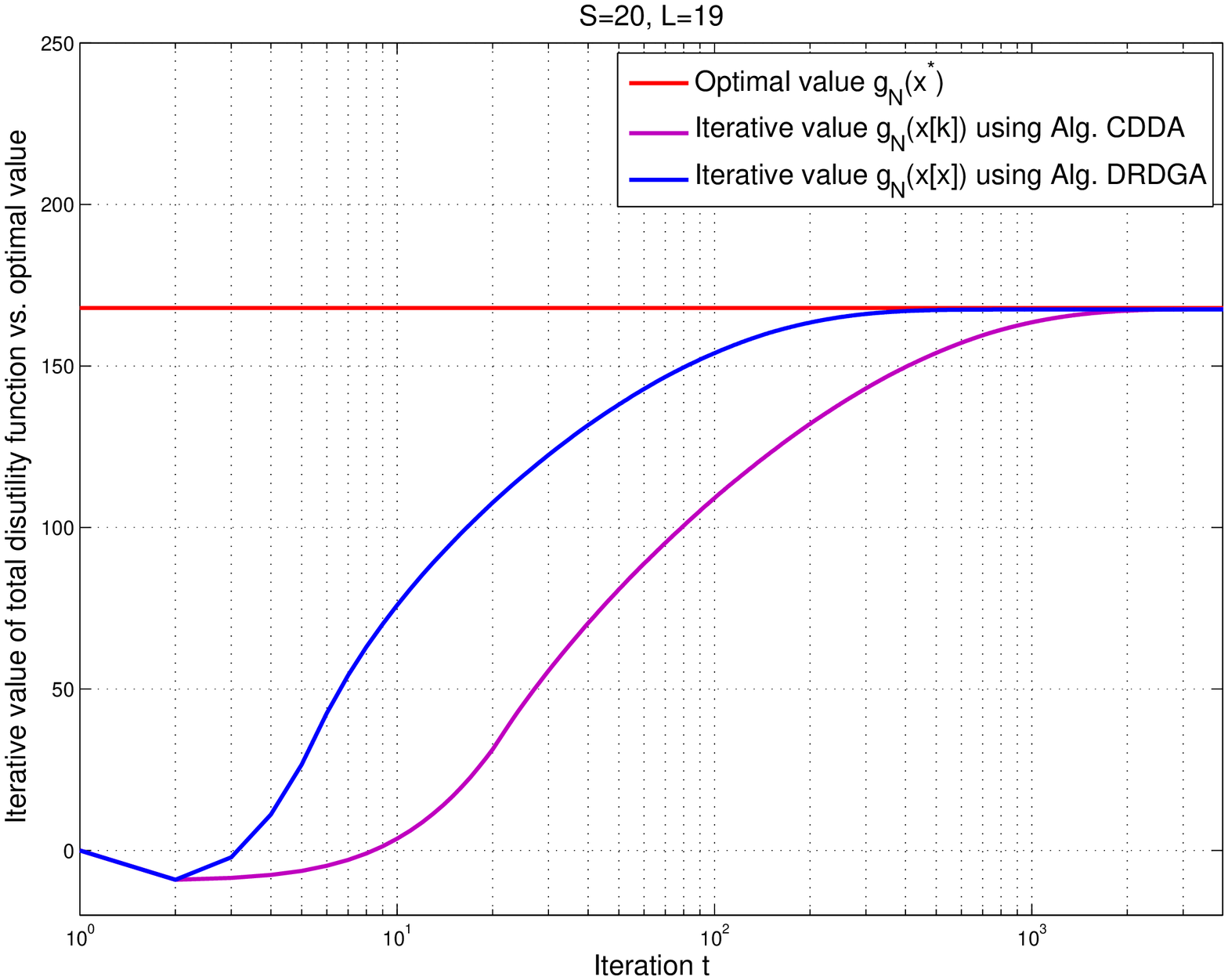}
\caption{Iterative value of total disutity function using Algs. DRDGA and CDDA vs. optimal value}
\label{fig6}
\end{figure}

\section{Conclusion}\label{6}

This paper proposed a solution tool for distributed convex problems with coupling equality constraints. The proposed algorithm is implemented in time-changing directed networks. By resorting to regularize the Lagrangian function, the norm of dual variables can be bounded. The proposed method can reach a fast convergence rate with order $O(\ln t/t)$ under some conditions. Numerical example on the network utility maximization demonstrates that the effectiveness of the proposed algorithm.
As a future research, it is interesting to analyze the communication delays of the proposed distributed method in this paper.


\section{References}\label{7}

\begin{thebibliography}{99}

\bibitem{nedic2009distributed}A. Nedi\'{c} and A. Ozdaglar. Distributed subgradient methods for multi-agent optimization. IEEE Transactions on Automatic Control, 2009,54(1):48-61.

\bibitem{nedic2010constaints}A. Nedi\'{c}, A. Ozdaglar, and P. Parrilo. Constrainted consensus and optimization in multi-agent networks. IEEE Transactions on Automatic Control, 2010,55(4):922-938.

\bibitem{jakovetic2014fast}D. Jakovetic, J. Xavier, and J.M. Moura. Fast distributed gradient methods. IEEE Transactions on Automatic Control, 2014,59(5):1131-1146.

\bibitem{nedic2015distributed}A. Nedi\'{c}, and A Olshevsky, Distributed optimization over time-varing directed graphs. IEEE Trans. Autom. Control, 2015,3(60):601-615.

\bibitem{johansson2008subgradient}B. Johansson, T. Keviczky, M. Johansson, and K. H. Johansson, Subgradient methods and consensus algorithms  for solving convex optimization problems, in Proc. IEEE CDC, Cancun, Mexico, Dec. 2008:4185-4190.

\bibitem{baingana2014proximalgradient}B. Baingana, G. Mateos, and G. Giannakis. Proximal-gradient algorithms for tracking cascades over social networks. IEEE Journal of Selected Topics in Signal Processing, 2014,8(4):563-575.

\bibitem{mateos2012distributed}G. Mateos and G. Giannakis. Distributed recursiveleast-squares: Stability and performance analysis. IEEE Transactions on Signal Processing, 2012,60(7):3740-3754.

\bibitem{bolognani2015distributed}S. Bolognani, R. Carli, G. Cavraro, and S. Zampieri. Distributed reactive power feedback control for voltage regulation  and loss minimization. IEEE Transactions on Automatic Control, 2015,60(4):966-981.

\bibitem{zhangdistributed2016}Y. Zhang and G. Giannakis. Distributed stochastic market clearing with high-penetration wind power and large-scale demand response. IEEE Transactions on Power Systems, 2016,31(2):895-906.

\bibitem{martinea2007on}S. Martinez, F. Bullo, J. Cortez, and E. Frazzoli. On synchronous robotic networks-Part I: Models, tasks, and complexity.  IEEE Transactions on Automatic Control, 2007,52(12):2199-2213.

\bibitem{tsitsiklis1986distributed}J. N. Tsitsiklis, D. P. Bertsekas, and M. Athans. Distributed asynchronous deterministic and stochastic gradient optimization algorithms. IEEE Transactions on Automatic Control, 1986,31(9):803-812.

\bibitem{ram2010distributed}S.S. Ram, A. Nedi\'{c}, and V.V. Veeravalli. Distributed stochastic subgradient projection algorithms for convex optimization. Journal of Optimization Theory and Applications, 2010,147(3):516-545.

\bibitem{c2012dual}J. C. Duchi, A. Agarwal, and M.J. Wainwright. Dual averaging for distributed optimization: convergence analysis and network scaling. IEEE Transactions on Automatic Control, 2012,57(3):592-606.

\bibitem{zhuon2012}M. Zhu, and S. Martinez. On distributed convex optimization under inequality and equality constraints. IEEE Transactions on Automatic Control, 2012,57(1):151-163.

\bibitem{li2015gradient}J. Li, C. Wu, Z. Wu, and Q. Long. Gradient-free method for nonsmooth distributed optimization. Journal of Global Optimization, 2015, 61(2):325-340.

\bibitem{lorenzo2016next}P.D Lorenzo, and G. Scutari. Netx: In-network nonconvex optimization. IEEE Transactions on Signal and Information Processing over Networks, 2016, 2(2):120-136.

\bibitem{gharesifard2012distributed}B. Gharesifard and J. Cortes. Distributed continuous-time convex optimization on weight-balanced digraphs. IEEE Transactions on Automatic Control, 2012,59(3):781-786.

\bibitem{iconsensus2012}K.I. Tsianos, S. Lawlor, and M.G. Rabbat. Consensus-based distributed optimization: Practical issues and applications in large-scale machine learning. In Communication, Control, and Computing (Allerton), 2012 50th Annual Allerton Conference on IEEE, 2012:1543-1550.

\bibitem{nedic2015stochastic}A. Nedi\'{c}, and A Olshevsky, Stochastic gradient-push for strongly convex functions on time-varying directed graphs. IEEE Transactions on Automatic Control, 2016,12(61):3936-3947.

\bibitem{bertsekas2003convex}D.P. Bertsekas, A. Nedi\'{c}, and A.E. Ozdaglar. Convex Analysis and Optimization. Belmont, MA, USA: Athena Scientific, 2003.

\bibitem{necoara2008application}I. Necoara, and J.A. Suykens. Application of smoothing technique to decomposition in convex optimization. IEEE Transactions on Automatic Control, 2008,53(11):2674-2679.

\bibitem{li2016a}J. Li, G. Chen, Z. Dong, and Z. Wu. A fast dual proximal-gradient method for separable convex optimization with linear coupled constraints. Computational Optimization and Applications, 2016,64(3):671-697.

\bibitem{yuan2011distributed}D. Yuan, S. Xu, and H. Zhao. Distributed primal-dual subgradient method for multiagent optimization via consensus algorithms. IEEE Transactions on Systems, Man, and Cybernetics, Part B (Cybernetics), ,2011,41(6):1715-1724.

\bibitem{chang2014distributed}T. H. Chang, A. Nedi\'{c}, and A. Scaglione. Distributed constrained optimization by consensus-based primal-dual perturbation method. IEEE Transactions on Automatic Control, 2014,59(6):1524-1538.

\bibitem{Yuan2016Regularized}D. Yuan, D.W.C. Ho, and S. Xu. Regularized primal-dual subgradient method for distributed constrained optimization. IEEE Transactions on Cybernetics, 46 (9): 2109-2118, 2016.

\bibitem{Khuzani2016Distributed}M.B. Khuzani, N. Li. Distributed regularized primal-dual method: Convergence analysis and trade-offs, arXiv preprint arXiv:1609.08262, 2016.

\bibitem{falsone2016dual}Alessandro. Falsone, Kostas. Margellos, Simone. Garetti, and Maria. Prandini. Dual decomposition and proximal minimization for multi-agent distributed optimization with coupling constraints. Automatica, 84 (2017) 149-158.

\bibitem{Low1999optimization}S. H. Low and D. E. Lapsley, Optimization flow control. I. basic algorithm and convergence. IEEE/ACM Transactions on Networking, 1999, 7:861-874.

\bibitem{beck2014an}A. Beck, A. Nedi\'{c}, A. Ozdaglar, and M. Teboulle. An $O({1}/{k})$ Gradient method for network resource Aalocation problems. IEEE Trans. Cont. Net. Sys, 2014,1(1):64-73.

\end{thebibliography}
\end{document}